\documentclass[10pt, final]{amsproc}
\usepackage[cp850]{inputenc}
\usepackage{amssymb}
\usepackage[mathscr]{eucal}
\usepackage{amscd}
\usepackage{amsmath}
\usepackage{setspace}
\usepackage{tikz}
\newcommand{\n}{nonnegative quadratic forms}

\def\Ker{\operatorname{Ker}}

\theoremstyle{plain}

\newtheorem{teor}{Theorem}
\newtheorem{prop}{Proposition}

\theoremstyle{remark}

\author{Jes\'us Su\'arez de la Fuente}
\address{Escuela Polit\'ecnica, Avenida de la Universidad s/n, 10071 C\'aceres, Spain.}
\email{jesus@unex.es}
\thanks{The author was supported in part by project MTM2016-76958-C2-1-P, project IB16056 and 
Ayuda a Grupos GR15152 de la Junta de Extremadura.}
\keywords{Grothendieck's theorem, quadratic form}
\subjclass{Primary: 46A22, 46B20}
\title{Subspaces of $\ell_1$ satisfying Grothendieck's theorem}

\begin{document} 

\begin{abstract} We characterize the subspaces $X$ of $\ell_1$ satisfying Grothendieck's theorem in terms of extension of nonnegative quadratic forms  $q:X \longrightarrow \mathbb R$ to the whole $\ell_1$.
\end{abstract}
\maketitle
\section{Introduction} The most famous contribution of Grothendieck to Banach space theory is perhaps what nowadays is called as ``Grothendieck's theorem". The survey written by Pisier \cite{Pi2} should convince any reader of the paramount importance of this result. It can be described as a nontrivial relation between the Hilbert space and the space $L_1$ (Theorem \ref{G}). Spaces verifying the theorem are called Grothendieck spaces, in short GT-spaces, that may be defined simply as those Banach spaces $X$ such that every operator $T:X\longrightarrow \ell_2$ is absolutely summing. Kalton and Pe\l czy\'nski characterized GT-spaces as follows (see \cite[Page 136]{KP} or \cite[Proposition 3.1.]{KP}):  Let $Z=Y/X$, where $Y$ is an $\mathcal L_1$-space, so that $X$ is the kernel of a quotient map onto $Z$. Then they showed that $X$ is a GT-space if and only if every short exact sequence 
$$
\begin{CD}
0 @>>>\ell_2@>>>\spadesuit@>>>Z@>>>0\;\;\;\;\\
\end{CD}
$$
splits. Much has been written about this kind of spaces although their structure is not still completely understood. As it is remarked in \cite[Remark 2, page 151]{KP}: \textit{It seems to be an interesting problem to characterize all Banach spaces $Z$ such that every twisted sum with $\ell_2$ as above splits.} It is not hard to see that this problem is equivalent to characterize the subspaces $X$ of $\ell_1$ for which every operator $X\longrightarrow  \ell_2$ is 1-summing that are precisely the subspaces of $\ell_1$ that are GT-spaces.

	Our purpose is to describe these subspaces of $\ell_1$ in terms of extension of nonnegative quadratic forms (Theorem \ref{main}).
\section{Background}\label{back}
\subsection{Quadratic forms}
All the Banach spaces considered throughout the paper will be real Banach spaces. Recall that given a real Banach space $X$, a function $q: X\longrightarrow  \mathbb R$ is a \textit{continuous quadratic form} if there exists a continuous bilinear form $b:X\times X\longrightarrow  \mathbb R$ such that $q(x)=b(x,x)$ for each $x\in X$. Recall also that an operator $T:X\longrightarrow  X^*$ is called \textit{symmetric} if $\langle Tx, y \rangle= \langle x, Ty \rangle$ for all $x,y\in X$. It is easy to see that the formula $$q(x)=\langle Tx, x\rangle $$ defines a one-to-one correspondence between continuous quadratic forms on $X$ and the symmetric operators. We are mainly interested in \n. Our main tool to study such functions is the following theorem of Kalton, Konyagin and Vesely (\cite[Theorem 1.2.]{KK}). Recall that a quadratic form is called \textit{delta-semidefinite} if it is the difference of two \n.
\begin{teor}\label{ds} Let $q$ be a quadratic form on a Banach space $X$ and $T:X\longrightarrow  X^*$ the symmetric operator that generates $q$. The following conditions are equivalent:
\begin{enumerate}
\item[(i)] $q$ is delta-semidefinite.
\item[(ii)] There exists a continuous quadratic form $p$ on $X$ such that $|q|\leq p$.
\item[(iii)] $T$ factors through a Hilbert space.
\end{enumerate}
\end{teor}
\subsection{The $p$-summing operators}
Let $T:X\longrightarrow Y$ be an operator between Banach spaces, and let $0<p< \infty$. Following \cite{Pi}, we will say $T$ is $p$-summing (in full, ``$p$-absolutely summing") if there is a constant $C$ such that, for all finite subsets $(x_i)\in X$, we have:
$$\left( \sum \|Tx_i\|^p  \right)^{\frac{1}{p}} \leq C \sup \left\{  \left (  \sum |\xi (x_i)|^p  \right)^{\frac{1}{p}} : \xi \in B(X^*) \right\}.$$
The least constant $C$ satisfying the inequality above is denoted by $\pi_p(T)$. Moreover, by $\Pi_p(X,Y)$ we denote the set of all $p$-summing operators. For $1\leq p <\infty$, $\Pi_p(X,Y)$ is a Banach space when endowed with the norm $\pi_p$. It will be useful for us to note that $\Pi_p(X,Y)$ satisfies the so-called ``\textit{ideal property}": If $T:X\to Y$ is $p$-summing and $\alpha:W\longrightarrow X$ and $\beta:Y\longrightarrow Z$ are bounded operators between Banach spaces, then the composition $\beta T \alpha$ is also $p$-summing and one has the estimate $$\pi_p(\beta T \alpha)\leq \|\beta \| \pi_p(T) \|\alpha\|.$$
We refer the reader to \cite{Pi} for further information on $p$-summing operators. Let us quote Grothendieck's theorem in the form it will be used in this paper. In what follows, let $H$ denote an arbitrary Hilbert space. 
\begin{teor}\label{G} Every bounded operator $T:L_1 \longrightarrow H$ is $1$-summing and satisfies $\pi_1(T)\leq K\|T\|$ for some absolute constant $K$.
\end{teor}
The theorem can be found in \cite[Theorem 5.10.]{Pi}. We will also use the following ``little Grothendieck's theorem":
\begin{teor}\label{g} Every bounded operator $T:C(K) \longrightarrow H$ is $2$-summing and satisfies $\pi_2(T)\leq \sqrt {\pi/2} \|T\|$.
\end{teor}
One may find this in \cite[Theorem 5.4.(b)]{Pi}. To finish this subsection recall again that GT-spaces are those Banach spaces $X$ such that every operator $T:X \longrightarrow H$ is $1$-summing. 
\subsection{Short exact sequences}
Recall that a \textit{short exact sequence} is a diagram like
\begin{equation}\label{ses}
\begin{CD} 0 @>>>X@>j>>Y@>\rho>>Z@>>>0
\end{CD}
\end{equation} where the morphisms are linear and such that the image of each arrow is the kernel of the next one. This condition implies that $X$ is a subspace of $Y$ and thanks to the open mapping theorem $Z$ is isomorphic to $Y/j(X)$. We say the sequence (\ref{ses}) \textit{splits} if there is a bounded linear map $R: Y\longrightarrow  X$ such that $R\circ j=Id_X$, where $R$ receives the name of a \textit{retraction} for $j$. Recall that in the case of bounded linear maps it is equivalent to have a retraction for $j$ than to have a selection for $\rho$. In this case, it is not hard to check that $Y$ is isomorphic to $X\oplus Z$.
\subsubsection{Push-out space, sequence and diagram}
 A commutative diagram
\begin{equation}\label{pushout}
\begin{CD}
X@>j>>Y\\
@ViVV @VI VV\\ 
V@>J>>PO,\\
\end{CD}
\end{equation}
is called a \textit{push-out} of 
$$
\begin{CD}
X@>j>>Y\\
@ViVV\\ 
V\\
\end{CD}
$$
provided that for every commutative diagram
$$
\begin{CD}
X@>j>>Y\\
@ViVV @V\beta VV\\ 
V@>\alpha>> W\\
\end{CD}
$$
there is a unique morphism $w:PO\longrightarrow  W$ so that $\alpha=wJ$, $\beta=w I$. There is a short description for the space $PO$ in (\ref{pushout}) called \textit{canonical push-out} and defined as $$PO=(Y\oplus_1 V)/\overline D$$
where $D=\{(jx,-ix):x\in X\}\subset
Y\oplus_1 V$. The $PO$ space is endowed with the natural quotient norm; and with $I$ and $J$ the compositions of the natural mappings of $Y$ and $V$ into $Y\oplus_1 V$ with the quotient map from $Y\oplus_1 V$ onto $PO$. Given the sequence (\ref{ses}) and the operator $i:X \longrightarrow
V$, we may complete (\ref{pushout}) to produce a commutative \textit{push-out diagram}
$$
\begin{CD}
0 @>>>X@>j>>Y@>\rho>>Z@>>>0\;\;\;\;&&(1)\\
&& @ViVV  @VIVV \|\\
0 @>>>V@>J>>PO@>P>>Z@>>>0\;\;\;\;&&\textit{(Push-out sequence).}\\
\end{CD}
$$
The sequence in the second row of the diagram above is called the \textit{push-out sequence} of the sequence (\ref{ses}) and $i$.  
\section{GT-spaces and nonnegative quadratic forms}\label{main1}
To prove our main result we need to characterize \n. Our equivalence follows directly from the proof of Theorem \ref{ds}.
\begin{prop}\label{non}
Let $q$ be a quadratic form on a Banach space $X$ and $T:X\longrightarrow  X^*$ the symmetric operator that generates $q$. The following conditions are equivalent:
\begin{enumerate}
\item[(i)] $q$ is a nonnegative quadratic form.
\item[(ii)] There exists a Hilbert space $H$ and an operator $A:X\longrightarrow  H$ such that $T=A^*A$.
\end{enumerate}
\end{prop}
\begin{proof}
(ii) $\implies$ (i) is obvious so let us prove the converse. If $q$ is nonnegative one may take $p=q$ in Theorem \ref{ds}. Therefore one find the following factorization for $T$: We first write $A:X\longrightarrow X/\Ker T$ as the natural quotient map where $X/\Ker T$ is the pre-Hilbert space (is not complete) given by the scalar product $$\langle [x], [y] \rangle= \langle Tx, y\rangle.$$ The operator $B: X/\Ker T \longrightarrow  X^*$ is given by the rule $B([x]):=Tx$. It is clearly well defined and the proof  that $B$ is bounded is similar to the one given in the proof of Theorem \ref{ds}. Therefore, completing  $X/\Ker T$ one has $T=BA$ (we denote in the same way the extension of $B$ to the completion) and thus 
$$q(x)=\langle BA x, x \rangle=\frac{1}{4}\| A(x)+B^*x\|^2-\frac{1}{4}\| A(x)-B^*x\|^2.$$
Observe that, in this particular version of Theorem \ref{ds}, one has $B^*_{|X}=A$ and then an easy computation finishes the proof.
\end{proof}
We are ready to prove our equivalence that is the main result of the paper.
\begin{teor}\label{main} Let $Z$ be a separable Banach space and let $\varphi:\ell_1 \longrightarrow  Z$ be any quotient map onto $Z$. Then the following are equivalent: 
\begin{enumerate}
\item[(i)] Every short exact sequence $0\to \ell_2 \to Y\to Z\to 0$ splits.
\item[(ii)] $\Ker \varphi$ is a GT-space.
\item[(iii)] Every nonnegative quadratic form  $q:\Ker \varphi\longrightarrow \mathbb R$ extends to a nonnegative quadratic form $Q: \ell_1\longrightarrow \mathbb R$.
\item[(iv)] Every nonnegative quadratic form $q:\Ker \varphi\longrightarrow \mathbb R$ extends to a delta-semidefinite quadratic form  $Q: \ell_1\longrightarrow \mathbb R$.
\end{enumerate}
\end{teor}
\begin{proof} The (i) $\iff$ (ii) is well known, see e.g. \cite[Proposition 3.1.]{KP}, but let us give a proof for the sake of completness. First of all, observe that given any short exact sequence \begin{equation}\label{one}
\begin{CD}
0 @>>>\ell_2@>j>>Y@>\rho>>Z@>>>0\;\;\;\;\\
\end{CD}
\end{equation}
one may construct the commutative diagram
\begin{equation}\label{po}
\begin{CD}
0 @>>>\Ker \varphi@>i>>\ell_1@>\varphi>>Z@>>>0\;\;\;\;\\
&& @VTVV  @VtVV \|\\
0 @>>>\ell_2@>j>>Y@>\rho>>Z@>>>0.\;\;\;\;\\
\end{CD}
\end{equation}
	It follows easily from the well known lifting property of $\ell_1$ that there is some $t$ making commutative the square on the right side and from this one may infer the existence of $T$. For specialists, (\ref{po}) is push-out diagram so some forthcoming arguments involving this fact and some extension and lifting principles may be simplified. However, our approach addresses also to a non specialist audience.\\
(i) $\implies $ (ii). Pick an arbitrary operator $T:\Ker \varphi \longrightarrow \ell_2$ and construct the push-out diagram 
$$
\begin{CD}
0 @>>>\Ker \varphi@>i>>\ell_1@>\varphi>>Z@>>>0\;\;\;\;\\
&& @VTVV  @VtVV \|\\
0 @>>>\ell_2@>j>>PO@>\rho>>Z@>>>0.\;\;\;\;\\
\end{CD}
$$
Since (i) holds there must a bounded retraction $r:PO \longrightarrow \ell_2$. Therefore $r\circ t$ is an extension of $T$ to $\ell_1$ and, by Theorem \ref{G}, the operator $r\circ t$ is $1$-summing. But, by the ideal property of $p$-summing operators, the operator $r\circ t \circ i=T$ is also $1$-summing.\\
(ii) $\implies$ (i). If $\Ker \varphi$ is a GT-space then, since 1-summing operators always have a bounded extension (they are also 2-summing \cite[Corollary1.6.]{Pi} and thus factor through an injective space, see e.g. Remark in \cite[Corollary1.8.]{Pi}), $T$ can be extended to the whole $\ell_1$ and this implies easily that the short exact sequence (\ref{one}) splits: If $\widehat{T}$ is such an extension, one readily verifies that $( j\circ \widehat{T} -t)\circ i=0$ and thus the operator  $j\circ\widehat{T}-t$ gives a well defined operator $S:Z \longrightarrow  Y$. It is a trivial computation to show that $\rho\circ S=Id_Z$ so that $S$ is a bounded linear selection for $\rho$.\\
(ii) $\implies $ (iii). Pick a nonnegative quadratic form $q: \Ker \varphi\longrightarrow  \mathbb R$ and use Proposition \ref{non} to show that the symmetric operator generating $q$ is of the form $A^*A$ for some $A:\Ker \varphi \longrightarrow  H$ and some Hilbert space $H$. Then one may construct the push-out diagram
$$
\begin{CD}
0 @>>>\Ker \varphi@>i>>\ell_1@>\varphi>>Z@>>>0\;\;\;\;\\\
&& @VAVV  @VaVV \|\\
0 @>>>H@>>>PO@>>>Z@>>>0.\;\;\;\;\\
\end{CD}
$$
Since by hypothesis the push-out sequence splits, the operator $A:\Ker \varphi \longrightarrow  H$ extends to an operator $\widehat{A}:\ell_1\longrightarrow  H$ (one just need, as before, to compose $a$ with a bounded linear retraction $r:PO \longrightarrow  H$). Then the extension of $q$ is $Q(x):=\|\widehat{A}(x)\|_{H}^2$.\\
(iii) $\implies$ (i). Pick any short exact sequence
$$
\begin{CD}
0 @>>>\ell_2@>j>>Y@>\rho>>Z@>>>0\;\;\;\;\\
\end{CD}
$$
and  construct again the commutative diagram (\ref{po})
$$
\begin{CD}
0 @>>>\Ker \varphi@>i>>\ell_1@>\varphi>>Z@>>>0\;\;\;\;\\
&& @VTVV  @VtVV \|\\
0 @>>>\ell_2@>j>>Y@>\rho>>Z@>>>0.\;\;\;\;\\
\end{CD}
$$
Let us define the quadratic form $$q(x)=\|T(x)\|^2,\;\;\;x\in \Ker \varphi.$$ Suppose $q$ extends to a nonnegative quadratic form $Q:\ell_1\longrightarrow  \mathbb {R}$. Then, by Proposition \ref{non}, $Q(x)=\|A(x)\|_{H}^2$ for some Hilbert space $H$ and certain operator $A:\ell_1\longrightarrow  H$.
Since $Q\circ i =q$ we find the equality $$\|T(x)\|_{\ell_2}=\|A\circ i(x)\|_{H},\;\;\;\;\;x\in \Ker \varphi.$$
Since $A$ is a 1-summing operator then so is $A\circ i$. Therefore, $T$ is also 1-summing and one may mimic the argument at the end of (ii) $\implies $(i). Observe that the equality above in particular gives $$\pi_p (T)=\pi_p(A\circ i),\;\;\;\;\;1\leq p<\infty.$$
(iii) $\implies$ (iv). Trivial.\\ 
(iv) $\implies$ (i). Pick a short exact sequence as in (\ref{one}) and  construct the push-out diagram (\ref{po}). Let us define the quadratic form $q(x)=\|T(x)\|^2$ and suppose $q$ extends to a delta-semidefinite quadratic form $Q:\ell_1\longrightarrow  \mathbb {R}$. Then by Theorem \ref{ds}, $T_Q$ -- the symmetric operator generated by $Q$ -- factors through some Hilbert space $H$ in the form $T_Q=BA$  in such a way that  $$Q(x)=\langle BA x, x \rangle=\frac{1}{4}\| A(x)+B^*x\|_H^2-\frac{1}{4}\| A(x)-B^*x\|_H^2.$$ 
Since $Q\circ i =q$ we obtain the equality $$\|T(x)\|_{\ell_2}^2=\frac{1}{4}\| A\circ i(x)+i^*\circ B^*x\|^2-\frac{1}{4}\| A\circ i (x)-i^*\circ B^*x\|^2,\;\;\;\;\;x\in \Ker \varphi.$$
Recall that $A$ and $B^*$ are 1-summing operators and then so are $A\circ i$ and $i^* \circ B^*$. This is enough to prove that $T$ must be also 1-summing:  Indeed, observe that the previous equality in particular gives
 \begin{eqnarray*}
  4\| T(x)\|^2 &\leq & \| A\circ i(x)+i^*\circ B^*x\|^2 + \| A\circ i (x)-i^*\circ B^*x\|^2\\
  &\leq & \left (\| A\circ i(x) \| +\| i^*\circ B^{*}x\| \right) ^2 +  \left ( \| A\circ i(x)\| +\| i^*\circ B^{*}x\|\right) ^2\\
  &= &  2\left ( \| A\circ i(x)\| +\| i^*\circ B^{*}x\|\right) ^2,
  \end{eqnarray*}
and thus
$$ \sqrt{2}\|T(x)\|\leq  \| A\circ i(x)\| +\| i^*\circ B^*x\|,$$
from where the claim follows. Actually one has $$\pi_1(T)\leq \frac{\pi_1(A)+\pi_1(B^*)}{\sqrt 2}.$$ Therefore $T$ can be extended and the sequence (\ref{one}) must split.
\end{proof}
The reader will have noticed that it is not really necessary to prove (iii) $\implies$ (i). However, we believe that, as the simplest case it is, it illustrates very well the key step of the proof.
\subsection{Remarks}
It is clear that interchanging the role of $\ell_1$ by an arbitrary GT-space $Y$, the same argument of Theorem \ref{main} gives a characterization of the subspaces of $Y$ that are GT-spaces.
Recall that Jarchow introduced \cite{J} the \textit{Hilbert-Schmidt spaces} as the spaces $X$ such that
$$\mathscr L (X,\ell_2)=\Pi_2(X, \ell_2).$$
\begin{prop}\label{hs} Let $X$ be a Banach space. Then the following conditions on $X$ are equivalent:
\begin{enumerate}
\item[(i)] $X$ is a Hilbert-Schmidt space.
\item[(ii)] Every nonnegative quadratic form on $X$ extends to a quadratic form on any space $Y$ containing $X$.
\end{enumerate}
\end{prop}
\begin{proof}
(i) $\implies$ (ii) is trivial so we prove the converse. Pick any  operator $T:X\longrightarrow \ell_2$ and consider the continuous quadratic form $q:X\longrightarrow \mathbb R$ given by $q(x)=\|T(x)\|^2$. Suppose $q$ extends to a quadratic form on $Y=C(B(X^*))$. Since every quadratic form on $C(K)$-spaces is delta-semidefinite \cite[Theorem 1.6]{KK}, a similar reasoning as in the proof of Theorem \ref{main} finishes the proof because $\mathscr L(C(K),\ell_2)=\Pi_2(C(K), \ell_2)$ by the little Grothendieck's theorem (Theorem \ref{g}).
\end{proof}


\begin{thebibliography}{99}
\bibitem{J} H. Jarchow, {\it On Hilbert-Schmidt spaces}, Proceedings of the 10th Winter School on Asbtract Analysis. Rend. Circ. Mat. Palermo \textbf{2} (1982), no. 2, 153--160.
\bibitem{KK} N. J. Kalton, S. V. Konyagin and L. Vesely, {\it Delta-semidefinite and delta-convex quadratic forms in Banach spaces}, Positivity \textbf{12} (2008), no. 2, 221--240.
\bibitem{KP} N. J. Kalton and A. Pe\l czy\'nski, {\it Kernels of surjections from $\mathcal L_1$-spaces with an application to Sidon sets}, Mathematische Annalen \textbf{309} (1997), no. 1, 135--158.
\bibitem{Pi} G. Pisier,  \emph{Factorization of linear operators and geometry of Banach spaces}, CBSM Regional Conference Series in Mathematics, 60.
\bibitem{Pi2} G. Pisier,  \emph{Grothendieck's theorem, past
and present}, Bull. Amer. Math. Soc. (N.S), 49 (2012), no. 2, 237--323.
\end{thebibliography}
\end{document}